\documentclass[aop,MSNbibl,citesort,dvips]{arximspdf}

\doi{10.1214/10-AOP599}
\volume{39}
\issue{2}
\pubyear{2011}
\firstpage{429}
\lastpage{434}

\makeatletter

\newtheorem{theorem}{Theorem}[section]

\makeatother

\begin{document}
\begin{frontmatter}

\title{T. E. Harris and branching processes}
\runtitle{T. E. Harris and branching processes}

\begin{aug}
\author[A]{\fnms{K. B.} \snm{Athreya}\corref{}\ead[label=e1]{kba@iastate.edu}
\ead[label=e2]{kbathreya@gmail.com}} and
\author[B]{\fnms{P. E.} \snm{Ney}\ead[label=e3]{ney@math.wisc.edu}}
\runauthor{K. B. Athreya and P. E. Ney}
\affiliation{Iowa State University and Indian Institute of Science, and
University of Wisconsin}
\address[A]{Departments of Mathematics and Statistics\\
Iowa State University\\
Ames, Iowa 50010\\
USA\\
\printead{e1}\\
and\\
Department of Mathematics\\
Indian Institute of Science\\
Bangalore, 560012\\
India\\
\printead{e2}}
\address[B]{Department of Mathematics\\
University of Wisconsin\\
Madison, Wisconsin 53705\\
USA\\
\printead{e3}}
\end{aug}

\received{\smonth{9} \syear{2010}}

%
\begin{abstract}
T. E. Harris was a pioneer par excellence in many fields of probability
theory. In this paper, we give a brief survey of the many fundamental
contributions of Harris to the theory of branching processes, starting
with his doctoral work at Princeton in the late forties and
culminating in his fundamental book ``The Theory of Branching
Processes,'' published in 1963.
\end{abstract}

%
\setattribute{keyword}{AMS}{AMS 2000 subject classification.}

\end{frontmatter}

\section{Introduction}

T. E. Harris wrote the first definitive
book \cite{H63} on branching processes, published in
1963. It covered much of the work on the subject up to that time,
a sizeable part due to Harris himself. It identified the
subject of branching processes and resulted in a great deal
of interest in the subject, among both mathematicians and
statisticians. Between
1963 and 1970, a vast number of papers on branching processes
appeared in many good journals specializing in
probability theory and mathematical statistics, and by
1971 more books on the subject appeared both in the U.S. and
elsewhere \cite{M71,S71}. Harris himself moved on to work on other beautiful
topics such as percolation and interacting particle systems. As with
branching processes, his work in these other areas was profound.
T. E. Harris
was pioneer par excellence, creating many areas of research in
which he laid the foundations that others built on. In what
follows, we present a brief account of Harris's contribution
to branching processes.

Harris's 1947 PhD dissertation at the mathematics department of Princeton
University was on branching processes,
titled ``Some theorems on Bernoulli multiplicative processes.''
This was followed in 1948 by his basic paper \cite{H48a} in the
\textit{Annals of
Mathematical Statistics}. In \cite{H48a}, he used the term
\textit{branching processes}, a~term which had also been used by
Russian mathematicians; he treated the single type discrete time
branching process. He also coined the term \textit{Galton--Watson branching
process} for this process. His main focus in \cite{H48a} was on the
supercritical case; we now give a description of this
work.

\section{Single type, discrete time case}

Let $\{p_j\}_{j\geq0}$ be a probability
distribution. Let $\{\xi_{n,k}; n\geq0,k\geq1\}$ be an array
of nonnegative integer valued random variables that are
i.i.d. (independent and identically distributed) with distribution
$\{p_j\}_{j\geq0}$. Let $Z_0$ be a positive integer. Now set
%
%
\begin{equation}
Z_1=\sum^{Z_0}_{k=1}\xi_{0,k}
\end{equation}
and for $n\geq1$,
$Z_{n+1}= \sum^{Z_n}_{k=1}\xi_{n,k}$ if $Z_n>0$ and 0
if $Z_n=0$.
Then the sequence $\{Z_n\}_{n\geq0}$ is called a Galton--Watson
branching process with initial population $Z_0$ and offspring
distribution $\{p_j\}_{j\geq0}$. Clearly, $\{Z_n\}_{n\geq0}$ is a
Markov chain with time homogeneous transition probabilities and
the nonnegative integers as the state space. The transition
probabilities are given by
\[
p_{ij}=P\Biggl(\sum^i_{r=1}\xi_r=j\Biggr) \qquad\mbox{for }
i\geq1\quad\mbox{and}\quad
p_{00}=1,
\]
where $\{\xi_r\}_{r\geq1}$ are i.i.d. with
distribution $\{p_k\}_{k\geq0}$.

One can interpret the sequence
$\{Z_n\}_{n\geq0}$ as follows. If $Z_n$ is thought of as the
number of individuals in the $n$th generation, then each one of them
produces a random number of children with distribution
$\{p_j\}_{j\geq0}$ independently of others in the $n$th generation
as well as any past ancestors. The total number $Z_{n+1}$ of
all these individuals is the size of the $(n+1)$st generation.

An important parameter in determining how the sequence
$\{Z_n\}_{n\geq0}$ behaves for $n$ large is the offspring mean
$m\equiv\sum_j j p_j$. Here are some basic results.
\begin{theorem}
Let $0<m\equiv\sum^\infty_{j=1} jp_j<\infty$ and let
$P(0<Z_0<\infty)=1$. Then:

\begin{enumerate}[(iii)]
\item[(i)] $m<1\Rightarrow P(Z_n\rightarrow
0$ as $n\rightarrow\infty)=1$,
\item[(ii)] $m=1, p_1<1\Rightarrow
P(Z_n\rightarrow0$ as $n\rightarrow\infty)=1$,
\item[(iii)] $m>1
\Rightarrow P(Z_n\rightarrow
0$ as $n\rightarrow\infty\mid Z_0=1)\equiv q<1$,
\end{enumerate}
where $q$ is the unique root of the equation
\[
s=f(s)\equiv\sum^\infty_{j=0} p_j s^j,\qquad 0\leq s<1.
\]
Further, $P(Z_n\rightarrow\infty$ as $n\rightarrow\infty\mid
Z_0=1)=1-q$, and
for any $k\geq1$, $P(Z_n\rightarrow
0$ as $n\rightarrow\infty\mid Z_0=k)=q^k$.
\end{theorem}

Harris in his book \cite{H63} notes in that in 1874 in \cite{WG74},
Galton and Watson did notice that
the \textit{extinction probability} $q$ satisfied $q=f(q)$, but
failed to notice that if $m>1$ the relevant root is less than one.
Galton and Watson's work was motivated by the problem of the survival
of British peerage names, posed by Galton in the \textit{London Times}
in the
1870s.

In his paper \cite{H48a}, which is based on his doctoral thesis, Harris
focused mainly on the
\textit{supercritical case}, that is, $m>1$. The case $m=1$ is
called \textit{critical} and $m<1$ is the \textit{subcritical case}. Let
$\{p_j\},m,\{Z_n\}$ be as in Theorem~1.
\begin{theorem}[(Supercritical case, \cite{H48a,H51})] Assume
$p_0=0,p_1<1,m>1$,
$ \sum^\infty_{j=1} j^2 p_j<\infty$ and
$0<Z_0<\infty$. Let $W_n\equiv Z_n/m^n,n\geq0$.
Then there exists a nonnegative random variable $W$ such that:
\begin{enumerate}[(iii)]
\item[(i)]$E((W_n-W)^2 \mid Z_0)\rightarrow0$ as $n\rightarrow\infty$,
\item[(ii)] $P(W=0)=0$,
\item[(iii)] $W$ has an absolutely continuous
distribution an $(0,\infty)$ with a continuous density,
\item[(iv)]
$E(W \mid Z_0=1)=1$.
\end{enumerate}
\end{theorem}

Harris \cite{H63} observes that J. L. Doob seems to have been the first to
note that $\{W_n\}_{n\geq0}$ is a martingale and, being
nonnegative, converges a.s. as $n\rightarrow\infty$.
Kesten and Stigum \cite{KS66} improved on this, as follows.
\begin{theorem}[\cite{KS66}] Let
$p_0=0,p_1<1,0<Z_0<\infty,1<m$, and $W_n=\frac{Z_n}{m^n}$.
Then:
\begin{longlist} [(ii)]
\item[(i)] $ \sum^\infty_1 j(\log
j)p_j<\infty\Rightarrow W_n\rightarrow W$ a.s. and in mean,
where $P(W=0)=0$, $E(W \mid Z_0=1)=1$ and $W$ has an absolutely continuous
distribution on $(0,\infty)$.

\item[(ii)] $ \sum^\infty_1 j(\log j)p_j=\infty\Rightarrow
W_n\rightarrow0$, a.s.
\end{longlist}
\end{theorem}

The work of A. N. Kolmogorov \cite{Ko38} in 1938 and A. M. Yaglom \cite
{Y47} in 1947 (see
\cite{H63}) led to the following.
\begin{theorem}[(Critical case)]
Suppose $m=1,p_1<1$ and $\sum^\infty_1 j^2 p_j<\infty$.
Then, as $n\rightarrow\infty$,
\begin{enumerate}[(ii)]
\item[(i)] $nP(Z_n>0 \mid Z_0=1)\rightarrow
\frac{\sigma^2}{2}$, where $\sigma^2\equiv\sum^\infty_1 j^2 p_j-1$,

\item[(ii)] $P(\frac{Z_n}{n}>x \mid Z_0=1,Z_n>0)\rightarrow
e^{-{2}/{\sigma^2}x}$, for all $0<x<\infty$.
\end{enumerate}
\end{theorem}
\begin{theorem}[(Subcritical case)] Let $m<1$. Then for all $j\geq1,\lim_n
P(Z_n=j \mid Z_n>0)\equiv b_j$ exists, $0<b_j<\infty$ and $\sum^\infty_{j=1}
b_j=1$.
\end{theorem}

In his book \cite{H63}, Harris presents extensions of Theorems 2.1,
2.2, 2.4 and 2.5 to the
multitype (finite type) case. In \cite{KS66}, Kesten and Stigum established
the analog of Theorem 2.3 above for the multitype Galton--Watson
process. See Athreya and Ney \cite{AN72} for details; see also
Sevastyanov \cite{S71} and Mode \cite{M71}.

\section{Single type, age dependent case}

In 1948 Harris, with Richard Bellman
\cite{H48b,H52}, formulated the theory of age dependent branching processes,
where each individual lives a random length of time and on death
creates a random number of individuals, and all individuals live
and reproduce independently of each other. Assuming all moments on
the offspring distribution and an absolutely continuous life time
distribution, they established an integral equation for the
probability generating function of $Z(t)$, the population size at
time $t$. They showed that in the supercritical case, $Z(t)e^{-\alpha t}$
converges in probability to a limit random variable $W$, where
$\alpha$ is the \textit{Malthusian parameter} defined by $m
\int^\infty_0 e^{-\alpha u}\,dG(u)=1$, with $G(\cdot)$ being the
distribution function of the lifetime of an individual. They
further showed that $W$ is nontrivial and has an absolutely continuous
distribution on $(0,\infty)$. There are analogs of Theorems 2.3 and
2.4 for this case, as well.

Conditions for the supercritical case were relaxed by later
authors; see
Athreya and Ney \cite{AN72}.

\section{General type case}

Harris also considered branching
processes with arbitrary type space by using the point process
approach. Here in any generation, one has a finite point process an
some type space $X$. The basic branching property of independent
production is retained. An individual located at $x \in X$ produces
children according to a point process over $X$ whose distribution
depends on $x$. All individuals act independently of each other. For
this, Harris used the method of moment generating functions. In
\cite{H60}, he established the analog of Theorem 2.2 in this context,
and applied this to nuclear cascades and related processes, as well as a
one-dimensional neutron model. For details on this, see Harris's
book \cite{H63}. Harris mentions that J. E. Moyal worked on similar
ideas. In the 1970s, Jagers and his colloborators in Sweden
developed this topic further in great detail (see \cite{Jag75}).
See also Ney \cite{N64a,N64b}.

\section{Cosmic-ray cascades}

Harris studied the theory of cosmic-rays cascades and supplemented
the work of nuclear physicists; Chapter 7 of his 1963 book \cite{H63}
deals with this topic. We present a brief summary of
Harris's work on cosmic-ray cascades as
discussed in his paper \cite{H57}. Here are the model assumptions:

(1) A photon of positive energy $\varepsilon$, moving through
homogeneous material, has probability $\lambda \,dt+o(dt)$ of being
transformed in the thickness interval $(t,t+dt)$
into two electrons, positive or negative, which receive energies
$\varepsilon U$ and $\varepsilon(1-U)$, respectively, where $U$ is a random
variable with an absolutely continuous distribution in $(0,1)$.
Note that the
role of time parameter is played by the thickness of the material.

(2) An electron loses (by ``collision'' or ``ionization'') a
deterministic amount of energy $\beta t$ in an interval of length
$t$.

(3) An electron radiates photons in such a way that the probability
that an electron of energy $\varepsilon$ emits a photon of energy
between $\varepsilon u$ and $\varepsilon(u+du)$ in a small thickness
interval of length $dt$ is $k(u)\,du\,dt$. Further, the energy that
goes to the radiated photon is subtracted from that of the parent
photon. A special case of interest for $k(\cdot)$ is
$k(u)=\frac{\mu}{u}+k_0(u)$ with $|\frac{dk_o(u)}{du}|\leq
c(1-u)^{-b},0<u<1$, where $c$ and $\mu$ are constants and $b<2$.

Here $\lambda$ and $\beta$ are constants independent of $t$ and
$\varepsilon$.

Under the above assumptions, Harris shows that if $\beta=0$ and
$\varepsilon_0(t)$ is the energy at time $t$ of an electron with
$\varepsilon_0(0)=1$ then, for $t>0,X_0(t)\equiv- \log\varepsilon_0(t)$
has an infinitely divisible distribution with probability density
$h_t(x),x>0$, whose characteristic function (Fourier transform) is
given by
\[
\int^\infty_0 e^{i\theta x}h_t(x)\,dx=\exp\biggl(t \int^\infty_0
(e^{i\theta u}-1)k(1-e^{-u})e^{-u}\,du\biggr).
\]
Harris\vspace*{-2pt} observes that the special case when
$k(u)=-\frac{\mu}{\log(1-u)}$ and $h_t(x)$ is the Gamma density
$\frac{x^{\mu t-1}e^{-x}}{\Gamma(\mu t)}$, this was given by Bhabha and
Heitler \cite{BH37}.\vspace*{1pt}

Next, Harris considers the random process $N(\varepsilon,t)$, the
total number of electrons at $t$ whose energies are greater than
$\varepsilon$, for $\varepsilon>0$.
Let $f_i(s,\varepsilon,t)\equiv E_i(s^{N(\varepsilon,t)}),0\leq s\leq
1$, where $E_i$ stands for expectation when the starting particle
is of energy 1, and is a photon for $i=1$ and an electron for
$i=2$. Harris shows that the following integro-differential
equation holds:
\[
\frac{\partial f_2}{\partial
t}(s,\varepsilon,t)=\int^1_0 \biggl[f_1\biggl(s,\frac{\varepsilon}{u},t\biggr)f_2
\biggl(s,\frac{\varepsilon}{1-u},t\biggr)-f_2(s,\varepsilon,t) \biggr]k(u)\,du,
\]
with $f_1(s,\varepsilon,0)=f_1(s,1,t)=f_2(s,l,t)=1$ for $t>0$ and
$f_2(s,\varepsilon,0)=s$ for \mbox{$\varepsilon<1$}.
Harris shows that the earlier results of Bartlett and Kendall \cite{BK51}
and of Janossy \cite{Jan50} could be deduced from the above.

Harris introduces a vector valued Markov process
$(I(t),\zeta(t)),t\geq0$, where $I(t)$ is the condition of a
single particle at time $t$ which can be a photon $(I=1)$ or an
electron $(I=2)$ and has energy $\zeta(t)$. He then derives the
limiting distribution of the process $(I(t),\zeta(t))$ as
$t\rightarrow\infty$ (assuming $\beta=0$) and is able to deduce
the earlier results of other authors as special cases.

Harris also obtained results for cascades with $\beta>0$. In
particular, he shows that when $\beta>0$, the energy $\varepsilon_1(t)$
of an
electron at time $t$ can be represented by
\[
\varepsilon_1(t)=\max\biggl\{0,\varepsilon_0(t)\biggl(1-\beta\int^t_0
\frac{ds}{\varepsilon_0(s)}\biggr)\biggr\}.
\]

\section{Concluding remarks}

T. E. Harris was deeply involved in the development of all aspects
of contemporary branching process theory. He laid a rigorous
foundation to areas where it had been lacking. His 1963 book \cite{H63}
is a beautiful and major work of scholarship. One can substantially credit
to its publication the explosion of work on branching
processes in the 1960s and 1970s and up to the present. It set
the impetus and direction of research on the subject for many
years. The present authors owe T. E. Harris a deep debt of gratitude
for this.

%

%
\printaddresses


\begin{thebibliography}{20}

\bibitem{AN72}
%
\begin{bbook}[vtex]
\bauthor{\bsnm{Athreya},~\bfnm{Krishna~B.}\binits{K.~B.}} \AND
\bauthor{\bsnm{Ney},~\bfnm{Peter~E.}\binits{P.~E.}}
(\byear{1972}).
\btitle{Branching Processes}.
\bpublisher{Springer}, \baddress{New York}.
\bid{mr={0373040}}
\end{bbook}
%
\endbibitem

\bibitem{BK51}
%
\begin{barticle}[mr]
\bauthor{\bsnm{Bartlett},~\bfnm{M.~S.}\binits{M.~S.}} \AND
\bauthor{\bsnm{Kendall},~\bfnm{David~G.}\binits{D.~G.}}
(\byear{1951}).
\btitle{On the use of the characteristic functional in the analysis of some
stochastic processes occurring in physics and biology}.
\bjournal{Proc. Cambridge Philos. Soc.}
\bvolume{47}
\bpages{65--76}.
\bid{mr={0039947}}
\end{barticle}
%
\endbibitem

\bibitem{H48b}
%
\begin{barticle}[vtex]
\bauthor{\bsnm{Bellman},~\bfnm{Richard}\binits{R.}} \AND
\bauthor{\bsnm{Harris},~\bfnm{Theodore~E.}\binits{T.~E.}}
(\byear{1948}).
\btitle{Age-dependent stochastic branching processes}.
\bjournal{Proc. Natl. Acad. Sci. USA}
\bvolume{34}
\bpages{601--604}.
\bid{mr={0027466}}
\end{barticle}
%
\endbibitem

\bibitem{BH37}
%
\begin{barticle}[auto:SpringerTagBib|2010-11-08|14:30:47]
\bauthor{\bsnm{Bhabha},~\bfnm{H.~J.}\binits{H.~J.}} \AND
\bauthor{\bsnm{Heitler},~\bfnm{W.}\binits{W.}}
(\byear{1937}).
\btitle{The passage of fast electrons and the theory of cosmic showers}.
\bjournal{Proc. Roy. Soc. London Ser. A}
\bvolume{159}
\bpages{432--458}.
\end{barticle}
%
\endbibitem

\bibitem{H63}
%
\begin{bbook}[vtex]
\bauthor{\bsnm{Harris},~\bfnm{Theodore~E.}\binits{T.~E.}}
(\byear{1963}).
\btitle{The Theory of Branching Processes}.
\bseries{Die Grundlehren der Mathematischen Wissenschaften}
\bvolume{119}.
\bpublisher{Springer}, \baddress{Berlin}.
\bid{mr={0163361}}
\end{bbook}
%
\endbibitem

\bibitem{H48a}
%
\begin{barticle}[vtex]
\bauthor{\bsnm{Harris},~\bfnm{T.~E.}\binits{T.~E.}}
(\byear{1948}).
\btitle{Branching processes}.
\bjournal{Ann. Math. Statist.}
\bvolume{19}
\bpages{474--494}.
\bid{mr={0027465}}
\end{barticle}
%
\endbibitem

\bibitem{H51}
%
\begin{binproceedings}[mr]
\bauthor{\bsnm{Harris},~\bfnm{T.~E.}\binits{T.~E.}}
(\byear{1951}).
\btitle{Some mathematical models for branching processes}.
In \bbooktitle{Proc. {S}econd {B}erkeley {S}ympos.
{M}ath. {S}tatist. {P}robab. 1950}
\bpages{305--328}.
\bpublisher{Univ. California Press}, \baddress{Berkeley, CA}.
\bid{mr={0045331}}
\end{binproceedings}
%
\endbibitem

\bibitem{H57}
%
\begin{barticle}[vtex]
\bauthor{\bsnm{Harris},~\bfnm{T.~E.}\binits{T.~E.}}
(\byear{1957}).
\btitle{The random functions of cosmic-ray cascades}.
\bjournal{Proc. Natl. Acad. Sci. USA}
\bvolume{43}
\bpages{509--512}.
\bid{mr={0087260}}
\end{barticle}
%
\endbibitem

\bibitem{H60}
%
\begin{bbook}[vtex]
\bauthor{\bsnm{Harris},~\bfnm{T.~E.}\binits{T.~E.}}
(\byear{1959}).
\btitle{On One-dimensional Neutron Multiplication. Research Memorandum RM-2317}.
\bpublisher{RAND Corporation}, \baddress{Santa Monica, CA}.
\end{bbook}
%
\endbibitem

\bibitem{H52}
%
\begin{barticle}[vtex]
\bauthor{\bsnm{Harris},~\bfnm{Theodore}\binits{T.}} \AND
\bauthor{\bsnm{Bellman},~\bfnm{Richard}\binits{R.}}
(\byear{1952}).
\btitle{On age-dependent binary branching processes}.
\bjournal{Ann. of Math. \textit{(2)}}
\bvolume{55}
\bpages{280--295}.
\bid{mr={0045971}}
\end{barticle}
%
\endbibitem

\bibitem{Jag75}
%
\begin{bbook}[vtex]
\bauthor{\bsnm{Jagers},~\bfnm{Peter}\binits{P.}}
(\byear{1975}).
\btitle{Branching Processes with Biological Applications}.
\bpublisher{Wiley}, \baddress{London}.
\bid{mr={0488341}}
\end{bbook}
%
\endbibitem

\bibitem{Jan50}
%
\begin{barticle}[mr]
\bauthor{\bsnm{J{\'a}nossy},~\bfnm{L.}\binits{L.}}
(\byear{1950}).
\btitle{On the absorption of a nucleon cascade}.
\bjournal{Proc. Roy. Irish Acad. Sect. A.}
\bvolume{53}
\bpages{181--188}.
\bid{mr={0045341}}
\end{barticle}
%
\endbibitem

\bibitem{KS66}
%
\begin{barticle}[mr]
\bauthor{\bsnm{Kesten},~\bfnm{H.}\binits{H.}} \AND
\bauthor{\bsnm{Stigum},~\bfnm{B.~P.}\binits{B.~P.}}
(\byear{1966}).
\btitle{A limit theorem for multidimensional {G}alton--{W}atson processes}.
\bjournal{Ann. Math. Statist.}
\bvolume{37}
\bpages{1211--1223}.
\bid{mr={0198552}}
\end{barticle}
%
\endbibitem

\bibitem{Ko38}
%
\begin{barticle}[auto:SpringerTagBib|2010-11-08|14:30:47]
\bauthor{\bsnm{Kolmogorov},~\bfnm{A.~N.}\binits{A.~N.}}
(\byear{1938}).
\btitle{On the solution of a biological problem}.
\bjournal{Tomsk Univ. Proc}
\bvolume{2}
\bpages{7--12}.
\end{barticle}
%
\endbibitem

\bibitem{M71}
%
\begin{bbook}[auto:SpringerTagBib|2010-11-08|14:30:47]
\bauthor{\bsnm{Mode},~\bfnm{C.~J.}\binits{C.~J.}}
(\byear{1971}).
\btitle{Multitype Branching Processes}.
\bpublisher{Elsevier}, \baddress{New York}.
\end{bbook}
%
\endbibitem

\bibitem{N64a}
%
\begin{barticle}[mr]
\bauthor{\bsnm{Ney},~\bfnm{P.~E.}\binits{P.~E.}}
(\byear{1964}).
\btitle{Generalized branching processes. {I}. {E}xistence and uniqueness
theorems}.
\bjournal{Illinois J. Math.}
\bvolume{8}
\bpages{316--331}.
\bid{mr={0162275}}
\end{barticle}
%
\endbibitem

\bibitem{N64b}
%
\begin{barticle}[mr]
\bauthor{\bsnm{Ney},~\bfnm{P.~E.}\binits{P.~E.}}
(\byear{1964}).
\btitle{Generalized branching processes. {II}. {A}symptotic theory}.
\bjournal{Illinois J. Math.}
\bvolume{8}
\bpages{332--350}.
\bid{mr={0162276}}
\end{barticle}
%
\endbibitem

\bibitem{S71}
%
\begin{bbook}[vtex]
\bauthor{\bsnm{Sevastyanov},~\bfnm{B.~A.}\binits{B.~A.}}
(\byear{1971}).
\btitle{Branching Processes}.
\bpublisher{Nauka}, \baddress{Moscow}.
\bid{mr={0345229}}
\end{bbook}
%
\endbibitem

\bibitem{WG74}
%
\begin{barticle}[vtex]
\bauthor{\bsnm{Watson},~\bfnm{H.~W.}\binits{H.~W.}} \AND
\bauthor{\bsnm{Galton},~\bfnm{F.}\binits{F.}}
(\byear{1874}).
\btitle{On the probability of extinction of families}.
\bjournal{J. Royal Anthropological Inst.}
\bvolume{6}
\bpages{138--144}.
\end{barticle}
%
\endbibitem

\bibitem{Y47}
%
\begin{barticle}[mr]
\bauthor{\bsnm{Yaglom},~\bfnm{A.~M.}\binits{A.~M.}}
(\byear{1947}).
\btitle{Certain limit theorems of the theory of branching random processes}.
\bjournal{Doklady Akad. Nauk SSSR (N.S.)}
\bvolume{56}
\bpages{795--798}.
\bid{mr={0022045}}
\end{barticle}
%
\endbibitem

\end{thebibliography}
\end{document}